\documentclass[11pt]{amsart}
\usepackage{amsmath}
\usepackage{cases}
\usepackage{mathrsfs}
\usepackage{amssymb}
\usepackage{amsfonts}
\usepackage{amssymb,amsfonts}

\allowdisplaybreaks

\begin{document}
\theoremstyle{plain}
\newtheorem{thm}{Theorem}[section]
\newtheorem{theorem}[thm]{Theorem}
\newtheorem{lemma}[thm]{Lemma}
\newtheorem{corollary}[thm]{Corollary}
\newtheorem{corollary*}[thm]{Corollary*}
\newtheorem{proposition}[thm]{Proposition}
\newtheorem{proposition*}[thm]{Proposition*}
\newtheorem{conjecture}[thm]{Conjecture}
%%%%%%%%%%%%%%%%%%%% Text roman %%%%%%%%%%%%%%%%%%%%%%%%%%%%%
\theoremstyle{definition}
\newtheorem{construction}{Construction}
\newtheorem{notations}[thm]{Notations}
\newtheorem{question}[thm]{Question}
\newtheorem{problem}[thm]{Problem}
\newtheorem{remark}[thm]{Remark}
\newtheorem{remarks}[thm]{Remarks}
\newtheorem{definition}[thm]{Definition}
\newtheorem{claim}[thm]{Claim}
\newtheorem{assumption}[thm]{Assumption}
\newtheorem{assumptions}[thm]{Assumptions}
\newtheorem{properties}[thm]{Properties}
\newtheorem{example}[thm]{Example}
\newtheorem{comments}[thm]{Comments}
\newtheorem{blank}[thm]{}
\newtheorem{observation}[thm]{Observation}
\newtheorem{defn-thm}[thm]{Definition-Theorem}

\newcommand{\sM}{{\mathcal M}}

%%%%%%%%%%%%%%%%%%%%%%%%%%%%%%%%%%%%%%%%%%%%%%%%%%%%%%%%%%%%%%

\title{The higher order terms in asymptotic expansion of color Jones
  polynomials}
        \author{Shengmao Zhu}
        \address{Department of Mathematics and Center of Mathematical Sciences, Zhejiang University, Hangzhou, Zhejiang 310027, China}
        \email{zhushengmao@gmail.com}


%% NB There should be only one primary classification, and zero or
%more secondary classifications.

\begin{abstract}
Color Jones polynomial is one of the most important quantum
invariants in knot theory. Finding the geometric information from
the color Jones polynomial is an interesting topic.  In this paper,
we study the general expansion of color Jones polynomial which
includes the volume conjecture expansion and the
Melvin-Morton-Rozansky (MMR) expansion as two special cases.
Following the recent works on $SL(2,\mathbb{C})$ Chern-Simons
theory, we present an algorithm to calculate the higher order terms
in general asymptotic expansion of color Jones polynomial from the
view of A-polynomial and noncommutative A-polynomial. Moreover, we
conjecture that the MMR expansion corresponding to the abelian
branch of A-polynomial. Lastly, we give some examples to illustrate
how to calculate the higher order terms. These results support our
conjecture.

\end{abstract}
\maketitle \tableofcontents
\section{Introduction}
Let $J_N(\mathcal{K};q)$ be the {\em normalized colored Jones
polynomial} of a knot $\mathcal{K}$ colored by the $N$-dimensional
irreducible representation of $SU(2)$. Thus,
$J_N(\text{unknot};q)=1$, $J_1(\mathcal{K};q)=1$ for all
$\mathcal{K}$ and $J_{2}(\mathcal{K};q)$ is the Jones polynomial of
$\mathcal{K}$. $J_N(\mathcal{K}; q)$ is an important quantum
invariant in knot theory. People want to find the geometric
information from $J_{N}(\mathcal{K};q)$.  More precisely, let
$q=e^\frac{2\pi i}{k}$ and consider the following limit,
\begin{align*}
k, \ N\rightarrow \infty, \quad u=\pi i\frac{N}{k}\quad
\text{fixed}.
\end{align*}
Now the question is what's the behavior of the limit
\begin{align}
\lim_{N\rightarrow \infty} J_{N}(\mathcal{K};e^{\frac{2u}{N}}).
\end{align}

The first progress in this direction is the volume conjecture. Let
us briefly review it. R.M. Kashaev defined a knot invariant
associated with the quantum dilogarithm and integer $N$, denoted by
$\langle \mathcal{K} \rangle_N$. He conjectured that for any
hyperbolic knot $\mathcal{K}$ \cite{Ka}, when $N\rightarrow \infty$,
\begin{align}
|\langle\mathcal{K}\rangle|_N\sim_{N\rightarrow \infty}
\exp\left(\frac{N}{2\pi}Vol(M_\mathcal{K})\right)
\end{align}
where $M_{\mathcal{K}}$ is equal to the knot complement
$S^3\setminus \mathcal{K}$. $Vol(M_\mathcal{K})$ is the hyperbolic
volume of $M_{\mathcal{K}}$. Then, in \cite{MM}, H. Murakami and J.
Murakami proved that for any knot $\mathcal{K}$,
\begin{align*}
\langle \mathcal{K} \rangle_N=J_{N}(\mathcal{K};e^{\frac{2\pi
i}{N}}).
\end{align*}
Moreover, they generalized the volume definition at right hand side
of $(2)$ to simplicial volume of any knot complement
$M_{\mathcal{K}}$. Now, the volume conjecture is formulated as
follow \cite{MM},
\begin{conjecture}[Volume conjecture]
For a knot $\mathcal{K}$,
\begin{align*}
|J_{N}(\mathcal{K};e^{\frac{2\pi i}{N}})|\sim_{N\rightarrow \infty}
\exp\left(\frac{N}{2\pi}Vol(M_\mathcal{K})\right)
\end{align*}
where $Vol(M_{\mathcal{K}})$ is the simplicial volume of knot
complement $M_{\mathcal{K}}=S^3\setminus \mathcal{K}$. In
particular, when $\mathcal{K}$ is hyperbolic, $Vol(M_{\mathcal{K}})$
is the hyperbolic volume of $M_{\mathcal{K}}$.
\end{conjecture}
We remark that the original volume conjecture was proposed for link
$\mathcal{L}$ \cite{Ka,MM}, but in this paper, we only consider the
case of knot $\mathcal{K}$. It is also possible to remove the
absolute value to consider the complexified volume conjecture
\cite{MMOT}: For any hyperbolic knot $\mathcal{K}$,
\begin{align}
J_{N}(\mathcal{K};e^{\frac{2\pi i}{N}})\sim_{N\rightarrow
\infty}\exp\left(\frac{N}{2\pi}\left(Vol(M_{\mathcal{K}})+iCS(M_{\mathcal{K}})\right)\right)
\end{align}
where $CS(M_{\mathcal{K}})$ is the Chern-Simons invariant of
$M_{\mathcal{K}}$ \cite{CS}. Furthermore , S. Gukov proposed a
$u$-parameterized version of complexified volume conjecture for any
hyperbolic knot $\mathcal{K}$ \cite{Gukov},
\begin{align}
J_{N}(\mathcal{K};e^{\frac{2u}{N}})\sim_{N\rightarrow
\infty}\exp\left(\frac{k}{\pi i}S_0(u)\right)
\end{align}
where $S_0(u)$ is a geometric invariant related the $u$-deformation
volume of $M_\mathcal{K}$ \cite{Thurston2}. In fact, formula $(4)$
is a generalization of $(3)$ for $u$ near the point $\pi i$ in
$\mathbb{C}$. Moreover, the expansion form of $(4)$ has been extend
to the higher order terms by S. Gukov and H. Murakami
\cite{Gukov-Murakami}.

It is also interesting to consider the situation when $u$ near the
point $0$. Another expansion form of color Jones polynomial called
the Melvin-Morton-Rozansky (MMR) conjecture was proposed in
\cite{MeM} and generalized by \cite{Rozansky}. The MMR conjecture
has been proved by D. Bar-Natan and S. Garoufalidis in \cite{BG}.
Recently, S. Garoufalidis and T. T. Q. Le obtained the following
analytic version of MMR expansion.
\begin{theorem}[\cite{GL2}]
For every knot $\mathcal{K}$, there exist a neighborhood
$\mathcal{O}\subset \mathbb{C}$ at $u=0$, such that for any $u\in
\mathcal{O}$, we have
\begin{align}
J_N(\mathcal{K};e^{\frac{2u}{N}})\sim_{N\rightarrow
\infty}\sum_{d=0}^{\infty}\frac{P_{\mathcal{K},d}(e^{2u})}{\Delta_\mathcal{K}(e^{2u})^{2d+1}}\left(\frac{2u}{N}\right)^d,
\end{align}
where $\Delta_\mathcal{K}(t)$ is the Alexander polynomial of
$\mathcal{K}$, $\{ P_{\mathcal{K},d}(t), d\geq 0\}$ is a sequence of
Laurent polynomials with $P_{\mathcal{K},0}(t)=1$.
\end{theorem}

Now we focus on the general expansion form of the limit $(1)$, it
was conjectured in \cite{Gukov} that the perturbative expansion of
$J_{N}(\mathcal{K};q)$ at the limit $N\rightarrow \infty, \
q\rightarrow 1$ was equal to the $SL(2,\mathbb{C})$ Chern-Simons
partition $Z(M_{\mathcal{K}})$ function up to a certain
normalization. Based on the standard perturbative Chern-Simons
theory \cite{AS1, Bar-Witten, AS2}, the general perturbative
computations of $Z(M_{\mathcal{K}})$ were explored in
\cite{DGLZ,DFM}. Therefore, motivated by the conjectured intimate
relation between the color Jones polynomial and Chern-Simons
partition, it is rational to consider the higher order expansion of
color Jones polynomial \cite{Gukov-Murakami}. If we introduce the
quantum parameter $\hbar$ as $ \hbar=\frac{i \pi}{k} $. The two
parameters $(k,N)$ in color Jones polynomial are changed to two
parameters $(\hbar, u)$. Then the general asymptotic expansion of
color Jones polynomial takes the following form \cite{DGLZ, DG},
\begin{align}
J_N(\mathcal{K};e^{\frac{2 u}{N}})\sim_{N\rightarrow
\infty}\exp\left(\frac{S_0(u)}{\hbar}-\frac{\delta_\mathcal{K}(u)}{2}\log\hbar+
\sum_{n=1}^{\infty}S_n(u)\hbar^{n-1}\right).
\end{align}

In this paper, we study the calculation of general terms $S_n(u)$
appearing above expansion $(5)$. We propose the two expansion
formulas $(4)$ and $(5)$ can be unified from the view of
$A$-polynomial and noncommutative $A$-polynomial of a knot
$\mathcal{K}$. In order to determine every terms $S_n(u)$ appearing
at the right side of $(6)$, one need to solve the following equation
with initial value $S_{\text{Initial}}(u)$:
\begin{equation} \label{eq:1}
\left\{ \begin{aligned}
&\hat{A}_{\mathcal{K}}(\hat{l},\hat{m};q)J_N(\mathcal{K};e^{\frac{2 u}{N}}) = 0 \\
                  &S_0(u)=S_{\text{Initial}}(u)
                          \end{aligned} \right.
                          \end{equation}
where the initial value $S_{\text{Initial}}(u)$ is determined by the
solution of the equation $A_{\mathcal{K}}(e^v,e^u)=0$ up to a
constant, where $A_{\mathcal{K}}(l,m)$ is the $A$-polynomial of
$\mathcal{K}$ and $\hat{A}_{\mathcal{K}}(\hat{l},\hat{m};q)$ is an
operator defined from the noncommutative $A$-polynomial
$\mathcal{K}$ which will be defined in section 2. More precisely, we
propose the following conjecture based on the work \cite{DGLZ}.
\begin{conjecture}
$i)$ There exists a solution of equation
$A_{\mathcal{K}}(e^v,e^u)=0$ called {\bf geometric branch} of
$A$-polynomial: $v=v^{\text{G}}(u)$. In this branch, we have a
neighborhood $\mathcal{O}^{G}\subset \mathbb{C}$ at $u=\pi i$, such
that for any $u\in \mathcal{O}^G$,
\begin{align}
J_N(\mathcal{K};e^{\frac{2 u}{N}})\sim_{N\rightarrow
\infty}\exp\left(\frac{S_0^{G}(u)}{\hbar}-\frac{3}{2}\log\hbar+
\sum_{n=1}^{\infty}S_n^{G}(u)\hbar^{n-1}\right),
\end{align}
with $\frac{dS_0^G(u)}{du}=v^{G}(u)$,  $S_0^{G}(u)$ is related to
the $u$-deformed volume of $M_\mathcal{K}$ by Gukov's conjecture
formula $(4)$, and
$S_1^{G}(u)=\frac{1}{2}\log\frac{iT_\mathcal{K}(u)}{4\pi}$
\cite{Gukov-Murakami}. Moreover, every $S_n^{G}(u)$ for $n\geq 2$
can be obtained by the algorithm introduced in section 2.

$ii)$ By the properties of $A$-polynomial, we know that there exist
an {\bf abelian branch} which corresponding to the branch $l=1$ of
$A_{\mathcal{K}}(l,m)=0$. In this branch, we have a neighborhood
$\mathcal{O}^A\subset \mathbb{C}$ of $0$, such that for any $u\in
\mathcal{O}^A$,
\begin{align}
J_N(\mathcal{K};e^{\frac{2 u}{N}})\sim_{N\rightarrow
\infty}\exp\left(\frac{S_0^{A}(u)}{\hbar}+
\sum_{n=1}^{\infty}S_n^{A}(u)\hbar^{n-1}\right),
\end{align}
with $S_0^A(u)=0$ and $S_1^A(u)=\log\frac{1}{
\Delta_\mathcal{K}(2u)}$, where $\Delta_{\mathcal{K}}(t)$ is the
Alexander polynomial of knot $\mathcal{K}$. Moreover, every
$S_n^{A}(u)$ for $n\geq 2$ can also be obtained by the same
algorithm.
\end{conjecture}

Furthermore, we promote that
\begin{conjecture}
The expansion formula $(9)$ is consistent with the analytic version
of MMR expansion $(5)$.
\end{conjecture}

\begin{remark}
In fact, there exists a sequence of Laurent polynomials $\{
Q_{\mathcal{K},n}(t)\}$ such that for $n\geq 2$, we have
\begin{align*}
S_n^{A}(u)=\frac{Q_{\mathcal{K}}(e^{2u})}{\Delta_{\mathcal{K}}(e^{2u})}.
\end{align*}
By the consistence of $(9)$ and  $(5)$,  if we let
$C_{\mathcal{K},d}(u)=\frac{2^d
P_{\mathcal{K},d}(e^{2u})}{\Delta_\mathcal{K}(e^{2u})^{2d+1}}$, then
\begin{align*}
C_{\mathcal{K},d}(u)=\exp (S_1^{A}(u))\sum_{\mu \mapsto
d}\frac{\prod_{i=1}^{l(\mu)}S_{\mu_i+1}^{A}(u)}{|Aut(\mu)|},
\end{align*}
where $\mu$ is the partition of $d$ with length $l(\mu)$. In other
words, conjecture 1.3 provides a method to compute every
$P_{\mathcal{K},d}(e^{2u})$ appears in the analytic version of MMR
expansion.
\end{remark}
\begin{remark}
The algorithm mentioned above to calculate the higher order terms
$S_n(u)$ is extracted from the recent works on perturbative
computation of $SL(2,\mathbb{C})$ Chern-Simons theory
\cite{DGLZ,DFM}. We note that, by the definitions in their works,
the color Jones polynomial $J_{N}(\mathcal{K};e^{\frac{2\pi i}{k}})$
and $Z(M_{\mathcal{K}};u,\hbar)$ are only difference with a
normalization
$\frac{q^{\frac{N}{2}}-q^{-\frac{N}{2}}}{q^{\frac{1}{2}}-q^{-\frac{1}{2}}}$.
Thus, we have the similar calculations for
$J_{N}(\mathcal{K};e^{\frac{2\pi i}{k}})$. Besides the geometric and
abelian branches, the authors also introduced the conjugate branch
of $A$-polynomial. See \cite{DGLZ} for more details.
\end{remark}
The rest of this paper is organized as follows: In section 2,
 we review the definitions of $A$-polynomial, non-commutative $A$-polynomial, AJ conjecture for
color Jones polynomial and their recent progresses. Then, we
illustrate the quantization algorithm to compute $S_n(u)$ which was
introduced in \cite{DGLZ} to study the pertubative computation of
$SL(2,\mathbb{C})$ Chern-Simons theory. In section 3, we give some
examples to illustrate the calculations of the higher order terms
with the quantization algorithm. More precisely, we have calculated
the following examples:

$i)$ Figure-8 knot $4_1$ in both geometric and abelian branches
which has been computed in \cite{DGLZ,DFM} with three different
methods under the  context of $SL(2,\mathbb{C})$ Chern-Simons
theory.

$ii)$ Twist knots $5_2$ and $6_1$ in abelian branch.

Lastly, the results on the abelian branch support our Conjecture
1.4.

\section{$A$-polynomial, noncommutative $A$-polynomial and the quantization algorithm}
\subsection{$A$-polynomial $A_\mathcal{K}(l,m)$ of a knot $\mathcal{K}$}
Let us start with the review of definition of $A$-polynomial of a
knot $\mathcal{K}$ in $S^3$ \cite{CCGD}. Denoted by
$R(M)=Hom(\pi_1(M),SL(2,\mathbb{C}))$ the set of all homomorphisms
$\rho$ from $\pi_1(M)$ to $SL(2,\mathbb{C})$ where $M=S^3\setminus
\mathcal{K}$. Let $R_{U}(M)$ be the subset of $R(M)$ consisting of a
representation $\rho$ such that $\rho(\mu)$ and $\rho(\lambda)$ are
upper triangular matrices for a fixed meridian $\mu$ and longitude
$\lambda$ of $\mathcal{K}$. Then one can define a projection
$\xi=(\xi_\lambda,\xi_\mu): R_U(M)\rightarrow \mathbb{C}^2$ by
$\xi(\rho)=(l,m)$ for $\rho\in R_{U}(M)$ with
\begin{equation*}
\rho(\lambda)=
\begin{pmatrix} l & * \\ 0 & l^{-1} \end{pmatrix} , \
\rho(\mu)=
\begin{pmatrix} m & * \\ 0 & m^{-1} \end{pmatrix}.
\end{equation*}
The Zariski closure of $\xi(R_{U}(M))$ is an algebraic variety in
$\mathbb{C}^2$ and each of its irreducible components is a curve,
which is defined by zeros of polynomial with integer coefficients in
$l$ and $m$. Then the product of those defining polynomials is
defined as the {\bf $A$-polynomial} of knot $\mathcal{K}$. Note that
the $A$-polynomial of $\mathcal{K}$ has a factor $l-1$, which
corresponds to to abelian representations which related to the
Alexander polynomial of $\mathcal{K}$. Thus someone define the
$A$-polynomial $A_{\mathcal{K}}(l,m)$ as the original $A$-polynomial
divided by $l-1$. The $A$-polynomial reflects the geometric
properties of the knot $\mathcal{K}$. More algebraic properties of
$A$-polynomial are listed in \cite{Gukov}.

Many $A$-polynomial of knot has been computed by now. Here we give
the $A$-polynomial of two types of knots. For a $(p,q)$-torus knot
$\mathcal{K}_{p,q}$, the $A$-polynomial is given by \cite{CCGD}:
\begin{align*}
A_{\mathcal{K}_{p,q}}(l,m)=1+lm^{pq}.
\end{align*}
Denote by $\mathcal{K}_{p}$ $p\in \mathbb{Z}$ the  $p$-twist knot,
its $A$-polynomial was computed in \cite{HS}.

When $p\neq -1,0,1,2$, $A_{\mathcal{K}_p}(l,m)$ is given recursively
by

\makeatletter
\let\@@@alph\@alph
\def\@alph#1{\ifcase#1\or \or $'$\or $''$\fi}\makeatother
\begin{subnumcases}
{A_{\mathcal{K}_p}(l,m)=}
cA_{\mathcal{K}_{p-1}}(l,m)-d A_{\mathcal{K}_{p-2}}(l,m), &$p> 0$, \label{eq:a1}
\\
cA_{\mathcal{K}_{p+1}}(l,m)-d A_{\mathcal{K}_{p+2}}(l,m), &$p<
0$.\label{eq:a2} \nonumber
\end{subnumcases}
\makeatletter\let\@alph\@@@alph\makeatother where
\begin{align*}
&c = -l + l^2 + 2lm^2 + m^4 + 2lm^4 + l^2m^4 + 2lm^6 + m^8 -lm^8, \\
&d=m^4(l + m^2)^4,
\end{align*}
and with the initial conditions
\begin{align*}
A_{\mathcal{K}_2}(l,m)&=-l^2 + l^3 + 2l^2m^2 + lm^4 + 2l^2m^4 -lm^6-
l^2m^8 \\&+2lm^{10} + l^2m^{10} + 2lm^{12} + m^{14} - lm^{14},\\
A_{\mathcal{K}_1}(l,m)&=l+m^6,\\
A_{\mathcal{K}_0}(l,m)&=1,\\
A_{\mathcal{K}_{-1}}(l,m)&=-l + lm^2 + m^4 + 2lm^4 + l^2m^4 +
lm^6-lm^8.
\end{align*}
For example, by the recursion $(10)$,
\begin{align*}
A_{\mathcal{K}_{-2}}(l,m)&=l^2 - l^3 - 3 l^2 m^2 + l^3 m^2 - 2 l m^4
- l^2 m^4 + 3 l m^6 + 3 l^2 m^6 \\&+ m^8 + 3 l m^8 + 6 l^2 m^8 + 3
l^3 m^8 + l^4 m^8 + 3 l^2 m^10 + 3 l^3 m^{10}\\& - l^2 m^{12} - 2
l^3 m^{12} + l m^{14} - 3 l^2 m^{14} - l m^{16} + l^2 m^{16}.
\end{align*}
\begin{remark}
The twist knots $\mathcal{K}_p$ for $p\in \mathbb{Z}$ include some
basic knots from Rolfsen's table.
\begin{align*}
&\mathcal{K}_1=3_1, \mathcal{K}_2=5_2, \mathcal{K}_3=7_2, \mathcal{K}_4=9_2, \\
&\mathcal{K}_{-1}=4_1, \mathcal{K}_{-2}=6_1, \mathcal{K}_{-3}=8_1,
\mathcal{K}_{-4}=10_1.
\end{align*}
\end{remark}
Recently, S. Garoufalidis and T. Mattman \cite{GaM} give a recursion
formula for the $A$-polynomial of the $(-2, 3, n)$ Pretzel knots.

\subsection{Noncommutative $A$-polynomial $\hat{A}_{\mathcal{K}}(E,Q;q)$}
The colored Jones polynomial $J_N(\mathcal{K};q)$ has many beautiful
structures. It was shown by S. Garoufalidis and TTQ Le \cite{Ga1,GL}
that the colored Jones function is $q$-holonomic, i.e. it satisfies
a nontrivial linear recursion relation with appropriate
coefficients. With such holonomicity, they introduce a geometric
invariant of a knot: the characteristic variety which is an affine
1-dimensional variety in $\mathbb{C}^2$. By comparing the character
variety of $SL(2,\mathbb{C})$ representations in the case of the
trefoil and figure-eight knots, they stated a conjecture that these
two varieties must be equal \cite{Ga1,GL}. They also define the
noncommutative $A$-polynomial $\hat{A}_{\mathcal{K}}(E,Q;q)$ for a
knot $\mathcal{K}$ which is the unique monic, linear, minimal order
$q$-difference equation satisfied by the sequence of color Jones
polynomials $\{J_{N}(\mathcal{K};q)\}$. Considering two operators
$E$ and $Q$ acting on the Jones polynomial $J_{N}(\mathcal{K};q)$ by
\begin{align}
(Q J_\mathcal{K})(N)=q^N J_{N}(\mathcal{K};q),
(EJ_{\mathcal{K}})(N)= J_{N+1}(\mathcal{K};q).
\end{align}
It is easy to see that $EQ=qQE$.

 Then $\hat{A}_{\mathcal{K}}(E,Q;q)$ controls
the recursion structure of color Jones polynomial
\begin{align}
\hat{A}_{\mathcal{K}}(E,Q;q)J_{N}(\mathcal{K};q)=0.
\end{align}
Note that $\hat{A}_{\mathcal{K}}(E,Q;q)$ can be written as the form
\begin{align}
\hat{A}_{\mathcal{K}}(E,Q;q)=\sum_{k\geq 0}a_k(Q;q) E^{k}
\end{align}
with $a_k(Q;q) \in \mathbb{Z}[q,Q]$. Then S. Garoufalidis
conjectured that
\begin{conjecture}[AJ Conjecture]
For every knot $\mathcal{K}$ in $S^3$,
$A_{\mathcal{K}}(l,m)=\epsilon \hat{A}_{\mathcal{K}}(l,m^2;q)$,
where $\epsilon$ is the evaluation map at $q=1$.
\end{conjecture}

In order to prove the AJ conjecture, a natural way is to compute the
non-commutative $A$-polynomial. So far, we have known an explicit
formula for torus knot in \cite{Ge}, figure-eight knot $4_1$ in
\cite{Ga1}, and 2-bridge knots \cite{Le}. Moreover, Takata found out
an explicit inhomogeneous $q$-difference equations for knots $5_2$
and $6_1$ with degree 5 and 6 respectively \cite{Ta}. But it is not
the really non-commutative $A$-polynomial in the sense of our
definition. Then S. Garoufalidis and X. Sun \cite{GS1,GS2} gave an
explicit formula for non-commutative $A$-polynomial of twist knots
$\mathcal{K}_p$ for $p=-8,..,11$. Recently, S. Garoufalidis and C.
Koutschan \cite{GaK} obtained the non-commutative $A$-polynomial the
for Pretzel knot $(-2,3,3+2p)$ for $p=-5,..,5$ using the method of
guessing.

Let us briefly describe the philosophy to calculate the
noncommutative $A$-polynomial $\hat{A}_{\mathcal{K}}(E,Q; q)$ of
knot $\mathcal{K}$.

For a generic planar projection of a knot $\mathcal{K}$, S.
Garoufalidis and T.T.Q. Le proved that the colored Jones polynomial
of a knot $\mathcal{K}$ can be written as a multisum \cite{GL}
\begin{align}
J_{N}(\mathcal{K};q)=\sum_{k_1,..,k_r}^{\infty}F(N,k_1,..,k_r),
\end{align}
where $F(N,k_1,..,k_r)$ is a proper $q$-hypergeometric function and
for a fixed $N\in \mathbb{Z}^+$, only finitely many terms are
nonzero. Because $F(N,k_1,..,k_r)$ is a proper $q$-hypergeometric
function, one can use the algorithm invented by Wilf-Zeilberger
\cite{PWZ,WiZeil}(the WZ algorithm ), also called creative
telescoping method,  to produce the noncommutative operator
eliminate $J_{N}(\mathcal{K};q)$. See \cite{PR,PR1} for a
mathematica implementation of $WZ$-algorithm. We will give some
examples to demonstrate how to use this computer program to derive
the noncommutative $A$-polynomial in next section.

\subsection{The algorithm to compute the asymptotic expansion of $J_{N}(\mathcal{K};q)$}
Let $A_{\mathcal{K}}(l,m)$ be the $A$-polynomial of a knot
$\mathcal{K}$. Define the operator $\hat{l}$ and $\hat{m}$ such that
\begin{align}
\hat{l}=E, \ \hat{m}^2=Q.
\end{align}
Then by $(12)$, we known that
$\hat{A}_{\mathcal{K}}(\hat{l},\hat{m}^2;q)J_{N}(\mathcal{K};q)=0$.

Recall the parameters we have described in the introduction section
\begin{align*}
\hbar=\frac{\pi i}{k}, \ u=\frac{\pi i N}{k}, \ q=e^{\frac{2\pi
i}{k}}.
\end{align*}
Then $q=e^{2\hbar}$,  the operator
$\hat{A}_{\mathcal{K}}(\hat{l},\hat{m}^2;q)$ annihilates
\begin{align*}
J(\mathcal{K};\hbar,u):=J_{N}(\mathcal{K};e^{2\hbar})
\end{align*}
i.e. we have the equation
\begin{align}
\hat{A}_{\mathcal{K}}(\hat{l},\hat{m}^2;q)J(\mathcal{K};\hbar,u)=0,
\end{align}
and by $(11)$ and $(15)$, the action of the operators $\hat{l},
\hat{m}$ is
\begin{align}
\hat{m}J(\mathcal{K};\hbar,u)=e^{u}J(\mathcal{K};\hbar,u),\quad
\hat{l}J(\mathcal{K};\hbar,u)=J(\mathcal{K};\hbar,u+\hbar).
\end{align}
It is clear that $\hat{l}\hat{m}=q^{\frac{1}{2}}\hat{m}\hat{l}$. As
in $(13)$, we expand $\hat{A}_{\mathcal{K}}(\hat{l},\hat{m}^2;q)$
as,
\begin{align}
\hat{A}_{\mathcal{K}}(\hat{l},\hat{m}^2;q)=\sum_{j=0}^{d}a_j(\hat{m},\hbar)\hat{l}^j.
\end{align}
Then we obtain
\begin{align}
\sum_{j=0}^{d}a_j(\hat{m},\hbar)J(\mathcal{K};\hbar,u+j\hbar)=0.
\end{align}
With the formula $(6)$, one can assume that at large $N$,
\begin{align*}
J(\mathcal{K};\hbar,u)=\exp\left(\frac{S_0(u)}{\hbar}-\frac{\delta_{\mathcal{K}}(u)}{2}\log\hbar+
\sum_{n=1}^{\infty}S_n(u)\hbar^{n-1}\right).
\end{align*}

 Therefore, from the above restriction equation for
 $J(\mathcal{K};q,u)$, one can obtain the sequence of expansion coefficients $\{
 S_{n}(u)\}$ recursively  by solving the equation $(16)$ for a given initial value $S_{0}(u)$ \cite{DGLZ}. In
following, we will show how to get the recursion formula for
$S_{n}(u)$ step by step.

Equation $(19)$ is equivalent to
\begin{align}
0&=\sum_{j=0}^{d}a_j(\hat{m},\hbar)\exp\left(\frac{1}{\hbar}S_{0}(u+j\hbar)-\frac{3}{2}\cdot
\log\hbar +\sum_{n\geq
0}\hbar^{n}S_{n+1}(u+j\hbar)\right)\\\nonumber
&=\exp\left(-\frac{\delta_{\mathcal{K}}(u)}{2}\cdot \log\hbar\right)
\sum_{j=0}^{d}a_j(\hat{m},\hbar)\exp\left(\sum_{n\geq -1}\hbar^n
S_{n+1}(u+j\hbar)\right).
\end{align}
And by Taylor expansion
\begin{align*}
\sum_{n\geq -1}\hbar^n S_{n+1}(u+j\hbar)&=\sum_{n\geq -1}\sum_{k\geq
0}\frac{S_{n+1}^{(k)}(u)}{k!}j^k\hbar^{k+n}\\\nonumber &=\sum_{t\geq
-1}\sum_{r=-1}^{t}\frac{S_{r+1}^{(t-r)}(u)}{(t-r)!}j^{t-r}\hbar^t\\\nonumber
&=\sum_{t\geq -1}S_{t+1}(u)\hbar^t+\sum_{t\geq
0}\sum_{r=-1}^{t-1}\frac{S_{r+1}^{(t-r)}(u)}{(t-r)!}j^{t-r}\hbar^{t}\\\nonumber
&=\sum_{t\geq -1}S_{t+1}(u)\hbar^t+jS_0'(u)+\sum_{t\geq
1}\sum_{r=-1}^{t-1}\frac{S_{r+1}^{(t-r)}(u)}{(t-r)!}j^{t-r}\hbar^{t}.
\end{align*}
It follows that
\begin{align*}
\sum_{j=0}^{d}a_j(\hat{m},\hbar)\exp\left(jS_{0}'(u)+\sum_{t\geq
1}B_t(u,j)\hbar^{t} \right)=0,
\end{align*}
where we have defined
\begin{align*}
B_t(u,j)=\sum_{r=-1}^{t-1}\frac{S^{(t-r)}_{r+1}(u)}{(t-r)!}j^{t-r}.
\end{align*}
Furthermore, one can expand $a_{j}(\hat{m},\hbar)$ and
$\exp\left(\sum_{t\geq 1}B_{t}(u,j)\hbar^{t}\right)$ as
\begin{align*}
a_{j}(\hat{m},\hbar)=\sum_{p\geq 0}a_{i,p}(\hat{m})\hbar^{p}
\end{align*}
and
\begin{align*}
\exp\left(\sum_{t\geq 1}B_{t}(u,j)\hbar^{t}\right)&=1+\sum_{n\geq
1}\frac{\left(\sum_{t\geq
1}B_{t}(u,j)\hbar^t\right)^n}{n!}\\\nonumber &=\sum_{\mu\in
\mathcal{P}}\frac{B_{\mu}(u,j)}{|Aut(\mu)|}\hbar^{|\mu|}
\end{align*}
where, $\mathcal{P}=\cup_{n\geq 1}\mathcal{P}_n\cup \{\emptyset \}$,
and $\mathcal{P}_n$ is the set of all partitions of integer $n\in
\mathbb{Z}^{+}$, and we denote $B_{\mu}(u,j)=B_{\mu_1}(u,j)\cdots
B_{\mu_l(\mu)}(u,j)$ and $B_{\emptyset }(u,j)=1$.

Then formula $(20)$ is equal to
\begin{align}
\sum_{j=0}^{d}e^{jS_{0}'(u)}\left(\sum_{p\geq 0}\sum_{\mu\in
\mathcal{P}}a_{j,p}(\hat{m})\frac{B_\mu(u,j)}{|Aut(\mu)|}\hbar^{|\mu|+p}\right)=0.
\end{align}
By the action of $\hat{m}$ defined in $(17)$, one can replace the
$\hat{m}$ with $e^u$ in $a_{j,p}(\hat{m})$.  As a series of $\hbar$,
all the coefficients of left hand side of $(21)$ must be zero.  The
constant term gives
\begin{align}
\sum_{j=0}^{d}e^{jS_{0}'(u)}a_{j,0}(e^u)=0.
\end{align}
which in fact is the $A$-polynomial.

When $n=|\mu|+p>0$, one can solve the $n$-th equation obtained from
the coefficient of $\hbar^{n}$ in equation $(21)$ and get
\begin{align}
S_n'(u)&=-\frac{1}{\sum_{j=0}^{d}e^{jS_{0}'(u)}a_{j,0}(e^u)j}\sum_{j=0}^{d}e^{jS_{0}'(u)}\left(
\sum_{p=1}^{n}a_{j,p}(e^u)\sum_{\mu\in \mathcal{P}_{n-p}\cup
\{\emptyset\}}\frac{B_{\mu}(u,j)}{|Aut(\mu)|}\right.\\\nonumber
&\left.+a_{j,0}(e^u)\sum_{\mu\in
\mathcal{P}\setminus\{(n)\}}\frac{B_{\mu}(u,j)}{|Aut(\mu)|}+a_{j,0}(e^u)\sum_{r=-1}^{n-2}\frac{S_{r+1}^{(n-r)}(u)}{(n-r)!}j^{n-r}
 \right).
\end{align}
\begin{example}
When $n=1$ and $2$, we have
\begin{align*}
S_1'(u)=-\frac{1}{\sum_{j=0}^{d}e^{jS_{0}'(u)}a_{j,0}(e^u)j}\sum_{j=0}^{d}e^{jS_{0}'(u)}\left(a_{j,1}(e^u)
+a_{j,0}(e^u)\frac{S_0''(u)}{2}j^2\right)
\end{align*}
\begin{align*}
S_2'(u)&=-\frac{1}{\sum_{j=0}^{d}e^{jS_{0}'(u)}a_{j,0}(e^u)j}\sum_{j=0}^{d}e^{jS_{0}'(u)}\left[a_{j,1}(e^u)
\left(\frac{S_0''(u)}{2}j^2+S'_1(u)j\right)\right.\\\nonumber
&\left.+a_{j,2}(e^u)
+a_{j,0}(e^u)\left(\frac{1}{2}\left(\frac{S''_0(u)}{2}j^2+S'_1(u)j\right)^2+
\frac{S'''_0(u)}{6}j^3+\frac{S_1''(u)}{2}j^2\right)\right]
\end{align*}
\end{example}

The above formula $(23)$ gives a recursion relation for $S_n(u)$. In
other words, if one knows the initial value $S_0(u)$, then all the
coefficients $S_{n}(u)$ are determined uniquely. How to choose
$S_0(u)$ depends on the choice of the branch of $A$-polynomial as
described in the introduction  section, i.e. in the geometric
branch: choosing $S_0(u)=S_0^G(u)$; and in abelian branch: choosing
$S_0(u)=S_0^{A}(u)=0$.
\begin{remark}
By AJ-conjecture, the classical limit $q\rightarrow 1$ of
noncommutative $A$-polynomial is the $A$-polynomial. Thus, the
noncommutative $A$-polynomial can be considered as the quantization
of $A$-polynomial. So the above method to compute $S_n(u)$ can be
called quantization algorithm. All the information of color Jones
polynomial are implied in a hierarchy of equations $(22), (23)$. The
first one equation $(22)$ is the $A$-polynomial if we let
$l=e^{S'_0(u)}$. $A$-polynomial reflects some geometric information
of the knot complement  $M_\mathcal{K}$. Finding the geometric
meaning of the generic equations $(23)$ will be interesting.
\end{remark}
\begin{remark}
The above quantization algorithm was introduce in \cite{DGLZ} to
study the  $SL(2,\mathbb{C})$ Chern-Simons partition function of
$M_{\mathcal{K}}$. They assumed that the color Jones polynomial
$J_{N}(\mathcal{K};e^{\frac{2\pi i}{k}})$ and Chern-Simons partition
$Z(M_{\mathcal{K}};u,\hbar)$ are only difference with a
normalization
$\frac{q^{\frac{N}{2}}-q^{-\frac{N}{2}}}{q^{\frac{1}{2}}-q^{-\frac{1}{2}}}$.
So in order to get the quantization operator
$\tilde{A}(\tilde{l},\tilde{m})$ of $Z(M_{\mathcal{K}};u;\hbar)$
such that $
\tilde{A}(\tilde{l},\tilde{m})Z(M_{\mathcal{K}};u;\hbar)=0$, one
only needs to do some modifications on operator
$\hat{A}_{\mathcal{K}}(\hat{l},\hat{m}^2;q)$. See \cite{DGLZ,DFM}
for detail discussion.
\end{remark}

\section{Examples}
\begin{example}

When $p=-1$, $\mathcal{K}_{-1}=4_1$, this example has been
calculated in \cite{DGLZ}. But we still recalculate here as an basic
example to illustrate the application of above algorithm.

{\bf Step 1}. Finding the noncommutative $A$-polynomial of $4_1$.

We download Paule and Riese's qZeil.m and qMultiSum.m package
\cite{PR} and run them in Mathematica 7.0.

In[1]:= $<$$<$  c:/qZeil.m

 q-Zeilberger Package by Axel Riese--@RISC Linz -V 2.42
(02/18/05)

In[2]:= $<$$<$ c:/qMultiSum.m

qMultiSum Package by Axel Riese--@RISC Linz -V 2.52 (30-Jul-2010)

In[3]:=summandfigure8 =
 $q^{n k} qfac[q^{-n - 1}, q^{-1}, k] qfac[q^{-n + 1}, q, k]$

Out[3]:={\em $q^{k n} \text{qPochhammer}[q^{-1 - n}, 1/q, k]
\text{qPochhammer}[q^{1 - n}, q, k]$}

In[4]:=qZeil[summandfigure8, {k, 0, Infinity}, n, 2]

Out[4]:=
\begin{align*}
\text{SUM[n]} &== \frac{q^{-1 - n}(q + q^n)(-q + q^{2 n})}{-1 + q^n}
- \frac{(1 - q^{-2 + n}) (1 - q^{-1 + 2 n}) \text{SUM[-2 + n]}}{(1 -
q^n) (1 - q^{-3 + 2 n})} +\\\nonumber
 &\frac{q^{-2 - 2 n} (1 - q^{-1 + n})^2 (1 + q^{-1 + n}) (q^4 +
q^{4 n}-q^{3 + n} - q^{1 + 2 n} - q^{3 + 2 n} - q^{1 + 3 n})
\text{SUM[-1 + n]}}{(1 - q^n) (1 - q^{-3 + 2 n})}
\end{align*}

This is a second-order inhomogeneous recursion relation, we convert
it into a third-order homogeneous recursion relation:

In[5]:=MakeHomRec [\%, SUM[n]];

Converting to forward shifts:

In[6]:=Rec41=ForwardShifts[\%]
\begin{align*}
\text{Out[6]}:=& q^{5 + n} (q - q^{3 + n}) (q^3 - q^{3 + n}) (q +
q^{3 + n}) (q - q^{6 + 2 n}) (q^3 - q^{6 + 2 n})
\text{SUM[n]}\\\nonumber & - q^{-5 - n} (q - q^{3 + n}) (q^2 - q^{3
+ n}) (q^2 + q^{3 + n}) (q - q^{6 + 2 n}) (q^3 - q^{6 + 2
n})\\\nonumber & \times (q^8 - 2 q^{9 + n} + q^{10 + n} - q^{9 + 2n}
+ q^{10 + 2 n} - q^{11 + 2 n} + q^{10 + 3 n} - 2 q^{11+ 3 n} + q^{12
+ 4 n})\\\nonumber &\times \text{SUM[1 + n]} + q^{-4 - n} (q - q^{3
+n})^2 (q + q^{3 + n}) (q^3 - q^{6 + 2 n}) (q^5 - q^{6 + 2 n}) (q^4
+ q^{5 + n}\\\nonumber & - 2 q^{6 + n} - q^{7 + 2 n} + q^{8 + 2 n} -
q^{9 + 2 n} - 2 q^{10 + 3 n} + q^{11 + 3 n} + q^{12 + 4 n})\text{
SUM[2 + n]} +q^{4 + n} \\\nonumber &\times(q - q^{3 + n}) (-1 + q^{3
+ n}) (q^2 + q^{3 + n}) (q^3 -q^{6 + 2 n}) (q^5 - q^{6 + 2 n})
\text{SUM[3 + n]} == 0
\end{align*}

Converting it to operator form:

In[7]:=F=$\text{ToqHyper[Rec41[[1]] - rec41[[2]]]} /.
\{\text{SUM[N]}\rightarrow 1, \text{SUM[N $q^c_.$]} :> X^c\} /.
\text{N} \rightarrow Q$
\begin{align*}
\text{Out[7]}:=&q^5 Q (q - q^3 Q) (q^3 - q^3 Q) (q + q^3 Q) (q - q^6
Q^2) (q^3 - q^6 Q^2) \\& - \frac{1}{q^5 Q}(q - q^3 Q) (q^2 - q^3 Q)
(q^2 + q^3 Q) (q - q^6 Q^2) (q^3 - q^6 Q^2)\\&\times (q^8 - 2q^9 Q +
q^{10} Q - q^9 Q^2 + q^{10} Q^2 - q^{11} Q^2 + q^{10} Q^3 - 2 q^{11}
Q^3 + q^{12} Q^4) X\\&
 + \frac{1}{q^4 Q}(q - q^3 Q)^2 (q + q^3 Q) (q^3 - q^6 Q^2) (q^5 -q^6 Q^2) (q^4 + q^5 Q - 2 q^6 Q\\&
- q^7 Q^2 + q^8 Q^2 - q^9 Q^2 - 2 q^{10} Q^3 + q^{11} Q^3 +
q^{12}Q^4) X^2 \\&+ q^4 Q (q - q^3 Q) (-1 + q^3 Q) (q^2 + q^3 Q)
(q^3-q^6 Q^2) (q^5 - q^6 Q^2) X^3
\end{align*}

Then $F$ is the non-commutative $A$-polynomial of $4_1$ if we
replace $X$ by $E$.

{\bf Step 2}. Finding the operator $
\hat{A}_{4_1}(\hat{l},\hat{m};q)=\sum_{j=0}^{d}a_j(\hat{m},\hbar)\hat{l}^j.
$

Substituting $Q$ and $X$ by $m^2$ and $l$ respectively in $F$, we
get
\begin{align*}
\hat{A}_{4_1}(\hat{l},\hat{m})=\sum_{j=0}^{3}a_j(\hat{m},\hbar)\hat{l}^j
\end{align*}
where
\begin{align*}
\hat{a}_0(\hat{m},q)&=\hat{m}^2q^5(q-\hat{m}^2q^3)(q^3 - \hat{m}^2
q^3) (q +\hat{m}^2 q^3) (q-\hat{m}^4 q^6) (q^3 -\hat{m}^4 q^6)
\\
\hat{a}_1(\hat{m},q)&=\frac{1}{\hat{m}^2 q^5}(q - \hat{m}^2 q^3)
(q^2 - \hat{m}^2 q^3) (q^2 + \hat{m}^2 q^3) (q -\hat{m}^4 q^6) (q^3 - \hat{m}^4 q^6) \\
 &\times(q^8 - 2 \hat{m}^2 q^9 - \hat{m}^4 q^9 + \hat{m}^2 q^{10} +\hat{m}^4 q^{10} + \hat{m}^6 q^{10} - \hat{m}^4 q^{11} - 2
\hat{m}^6 q^{11} + \hat{m}^8 q^{12})
\\
\hat{a}_2(\hat{m},q)&=\frac{1}{\hat{m}^2 q^4}(q - \hat{m}^2 q^3)^2
(q + \hat{m}^2 q^3) (q^3 - \hat{m}^4 q^6) (q^5 -\hat{m}^4 q^6)\\
& \times(q^4 + \hat{m}^2 q^5 - 2 \hat{m}^2 q^6 - \hat{m}^4 q^7 +
\hat{m}^4 q^8 -\hat{m}^4 q^9 - 2 \hat{m}^6 q^{10} + \hat{m}^6 q^{11}
+ \hat{m}^8 q^{12})
\\
\hat{a}_3(\hat{m},q)&=\hat{m}^2q^4 (q - \hat{m}^2 q^3) (-1 +
\hat{m}^2 q^3) (q^2 + \hat{m}^2 q^3) (q^3 - \hat{m}^4 q^6) (q^5 -
\hat{m}^4 q^6)
\end{align*}

{\bf Step 3} Choosing the different branches.

The $A$-polynomial of $4_1$ is
\begin{align*}
A_{4_1}(l,m)=(-1 + l)(l - l m^2 - m^4 - 2 l m^4 -l^2 m^4 - l m^6 + l
m^8)
\end{align*}

Solving this equation, we obtain the three branches: $l_{A}=1$ is
called the abelian branch and $l_{G}=-\frac{-1 + m^2 + 2 m^4 + m^6 -
m^8 + (-1 + m^4) \sqrt{1 - 2 m^2 - m^4 - 2 m^6 + m^8}}{2 m^4}$ is
the geometric branch. The third one is the conjugate of $l_{G}$
called conjugate branch which have the intimate relation with
geometric branch discussed in \cite{DGLZ}.

{\bf Step 4} Calculating the expansion coefficients $S_{n}(u)$ in
different branches by formula $(23)$.

Abelian branch expansion: taking the initial value $S^{A}_0(u)=\log
l_A=0$, then
\begin{align*}
&S_1^{A}(u)=\log\frac{1}{\Delta_{4_1}(m^2)};  \\
&S_2^{A}(u)=\text{constant};\\
&S_3^{A}(u)=\frac{4 (m^{-2} - 1 + m^2)}{\Delta_{4_1}(m^2)^3};\\
&S_4^{A}(u)=\text{constant};\\
&\cdots
\end{align*}
where $\Delta_{4_1}(t)=\frac{1}{t}+t-3$ is the Alexander polynomial
of $4_1$. The above results match the Conjecture 1.4.

Geometric branch expansion: the initial value
$S^{G}_0(u)=\frac{i}{2}Vol(4_1)+\int_{i\pi}^{u}v_{G}(u)du-2\pi^2$
\cite{DGLZ}.
\begin{align*}
&S_1^{G}(u)=2 \log(m) - \log(-1 + m^2) - \frac{1}{4} \log(1 - 2 m^2 - m^4 - 2 m^6 + m^8); \\
&S_2^{G}(u)=\frac{1 - m^2 - 2 m^4 + 15 m^6 - 2 m^8 - m^{10} + m^{12}}{12(-1 + 2 m^2 + m^4 + 2 m^6 - m^8)^{\frac{3}{2}}};\\
&S_3^{G}(u)=-\frac{2 m^6 (-1 + m^2 + 2 m^4 - 5 m^6 + 2 m^8 + m^{10}-m^{12})}{(1 -2 m^2 - m^4 - 2 m^6 + m^8)^3};\\
&S_4^{G}(u)=\frac{m^2}{90(1 - 2 m^2 - m^4 - 2 m^6
+m^8)^\frac{9}{2}}(1 - 4 m^2 - 128 m^4 + 36 m^6 \\&+ 1074 m^8 - 5630
m^{10}+ 5782 m^{12} + 7484 m^{14} - 18311 m^{16} + 7484 m^{18}\\& +
5782 m^{20} - 5630 m^{22} +1074 m^{24} + 36 m^{26} - 128 m^{28} - 4 m^{30} + m^{32});\\
&\cdots
\end{align*}
If we use the Ray-Singer torsion of $4_1$ \cite{Gukov-Murakami,DG}
\begin{align*}
T_{4_1}(u)=\frac{4\pi^2 m^2}{\sqrt{-1+2m^2+m^4+2m^6-m^8}},
\end{align*}
we may conjecture that $S_n(u)$ has the form
\begin{align*}
S_n(u)=\left(\frac{T_{4_1}(u)}{4\pi^2}\right)^{3n-3}G_n(m)  \
\text{for}\ n\geq 2,
\end{align*}
where $\{ G_n(m)\}$ is a sequence of Laurent polynomial of $m$.

\begin{remark}
In \cite{DGLZ}, they have calculated the perturbative expansion for
$Z\left(M_{4_1};u;\hbar\right)$ assume that
$Z(M_{4_1};u;\hbar)=\frac{q^{N/2}-q^{-N/2}}{q^{1/2}-q^{-1/2}}J(\mathcal{K};u,\hbar)$.
By this relation, we should make a modification for
$\hat{a}_{j}(\hat{m},q)$,
\begin{align*}
\hat{a}_{j}(\hat{m},q)\rightarrow\frac{\hat{a}_{j}(\hat{m},q)}{m^2
q^{j/2}-q^{-j/2}}.
\end{align*}
With these new $\hat{a}(\hat{m},q)$, they calculated the $S_n(u)$
for $Z(M_{4_1};u,\hbar)$ up to $n=8$.
\end{remark}
\end{example}

We give more examples in abelian branch expansion.
\begin{example}
When $p=2$, the twist knot $\mathcal{K}_{2}=5_2$. Setting the
initial value $S_0^A(u)=0$, we get
\begin{align*}
&S_1^A(u)=\log\left(\frac{1}{\Delta_{5_2}(m^2)}\right);\\
&S_2^A(u)=\frac{-4 (m^{-2}+ m^2)+13}{2\Delta_{5_2}(m^2)^2};\\
&S_3^A(u)=-\frac{-32 + 104 m^2 + 200 m^4 - 607 m^6 + 200 m^8 +
104 m^{10} - 32 m^{12}}{8 m^6\Delta_{5_2}(m^2)^4};\\
&S_4(u)=-\frac{1}{24 m^{10}\Delta_{5_2}(m^2)^6}
 (320 - 752 m^2 - 3808 m^4 + 3052 m^6 + 39692 m^8 \\
 &- 78163 m^{10} +
    39692 m^{12}+ 3052 m^{14} - 3808 m^{16} - 752 m^{18} + 320
    m^{20})\\&
\cdots
\end{align*}
where $\Delta_{5_2}(t)=2(t^{-1}+t)-3$ is the Alexander polynomial of
$5_2$. It is easy to see these results match the Conjecture 1.4.
\end{example}
\begin{example}
When $p=-2$, the twist knot $\mathcal{K}_{-2}=6_1$. Setting the
initial value $S_0^A(u)=0$, we obtain
\begin{align*}
&S_1^{A}(u)=\log\left(\frac{1}{\Delta_{6_1}(m^2)}\right);\\
&S_2^{A}(u)=\frac{-4 m^2 + 7 m^4 - 4 m^6}{2m^4 \Delta_{6_1}(m^2)^2};\\
&S_3^{A}(u)=-\frac{32 - 504 m^2 + 1656 m^4 - 2303 m^6 + 1656 m^8 -
504m^{10} +32 m^{12}}{8 m^6\Delta_{6_1}(m^2)^4};\\
&S_4(u)=-\frac{1}{24m^{10} \Delta_{6_1}(m^2)^6}(320 - 2512 m^2 +
23968 m^4 - 103404 m^6 + 225900 m^8\\& - 288925 m^{10} + 225900
m^{12} -103404 m^{14} + 23968 m^{16} - 2512 m^{18} + 320 m^{20})\\
&\cdots
\end{align*}
where $\Delta_{6_1}(t)=2(t^{-1}+t)-5$ is the Alexander polynomial of
$6_1$. These results match the Conjecture 1.4.
\end{example}
\begin{remark}
In the above two examples, we have used the non-commutative
$A$-polynomial for twist knot $\mathcal{K}_p$ obtained by S.
Garoufalidis and X. Sun \cite{GS1,GS2}. In fact, they have
calculated all the non-commutative $A$-polynomial of $\mathcal{K}_p$
for $p=-8,..,11$. So we can get more examples by using their
results.
\end{remark}

\section{Conclusion and discussion}
In this paper, we present an algorithm to calculate the higher order
terms in general expansion of color Jones polynomial from the view
of A-polynomial and noncommutative A-polynomial. In the large $N$
limit, the color Jones polynomial $J_{N}(\mathcal{K};e^{2u})$ has
the following expansion form
\begin{align}
J(\mathcal{K};\hbar,u)=\exp\left(\frac{S_0(u)}{\hbar}-\frac{3}{2}\log\hbar+
\sum_{n=1}^{\infty}S_n(u)\hbar^{n-1}\right).
\end{align}
In order to determine every terms $S_n(u)$ appearing at left side of
$(24)$, we need to solve the following equation with initial value
$S_{\text{Initial}}(u)$:
\begin{equation} \label{eq:1}
\left\{ \begin{aligned}
&\hat{A}_{\mathcal{K}}(\hat{l},\hat{m};q)J(\mathcal{K};\hbar,u) = 0 \\
                  &S_0(u)=S_{\text{Initial}}(u)
                          \end{aligned} \right.
                          \end{equation}
Up to a constant, the initial value $S_{\text{Initial}}(u)$ is
determined by the solution of the equation
$A_{\mathcal{K}}(e^v,e^u)=0$, where $A_{\mathcal{K}}(l,m)$ is the
$A$-polynomial of knot $\mathcal{K}$. More precisely, if we assume
$A_{\mathcal{K}}(l,m)=(l-1)f_{d}(l,m)$, where
$f_d(l,m)=\sum_{i=1}^{d}a_i(m)l^i$. For a given $m$, the equation
$f_d(l,m)=0$ has $d$ solutions in $\mathbb{C}$ denoted by
$l=l^{\alpha}(m)$, $\alpha=1,..,d$. Thus $A_{\mathcal{K}}(l,m)=0$
has $d+1$ branches: abelian branch $l^A=1$, and $l^{\alpha}(m)$, for
$\alpha=1,..,d$. There are some symmetries between these different
branches $\alpha=1,..,d$. See \cite{DGLZ} for discussions from
Chern-Simon theory.

One of the most interesting branch is called geometric branch
denoted by $l^G(m)$ which is relevant with the hyperbolic volume of
knot complement $M_{\mathcal{K}}$. In this geometric branch, the
initial value is $\frac{dS_{\text{Initial}}(u)}{du}=v^{G}(u)$ and
$S_{\text{Initial}}(u)$ is the complexified volume of
$M_\mathcal{K}$ parameterized by $u$. Then all the terms $S_n^G(u)$
can be solved from the recursive relation $(23)$. Moreover,
$S_1^{G}(u)$ has the geometric interpretation
$S_1^{G}(u)=\frac{1}{2}\log\frac{iT_\mathcal{K}(u)}{4\pi}$
\cite{Gukov-Murakami}, where $T_\mathcal{K}(u)$ $u$-deformed torsion
of $M_\mathcal{K}$. But what's the geometric mean of $S_n^G(u)$ for
$n\geq 2$ is still unknown.

In the abelian branch $l_A=1$, the initial value is
$\frac{dS_{\text{Initial}}(u)}{du}=v^A=0$. So
$S_{\text{Initial}}(u)$ is a constant.  One can get every $S_n^A(u)$
by formula $(23)$. Moreover, the first term $S_1^A(u)=\log\frac{1}{
\Delta_\mathcal{K}(2u)}$, where $ \Delta_\mathcal{K}(t)$ is the
Alexander polynomial of $\mathcal{K}$ which has the geometric
meaning, but the geometric interpretation is still unclear for
$S^{A}_n(u),n\geq 2$. We found that expansion terms $S_n^A(u)$ is in
consistence with the Melvin-Morton-Rozansky expansion for color
Jones polynomial \cite{MeM,DG,Rozansky,Rozansky2,GL}, so we proposed
the Conjecture 1.4.

We have calculated the following examples in this paper: i) figure-8
knot $4_1$ in both geometric and abelian branches; ii) Twist knots
$5_2$ and $6_1$ in abelian branch. Both of them match the
conjecture. For the cases in ii) , we only compute the abelian
branch, because the
 geometric branch in this situation is more complicated which consumes much computer time.

The recursive algorithm to get $S_n(u)$ was found in \cite{DGLZ} to
study the $SL(2,\mathbb{C})$ Chern-Simons partition function
$Z(M_{\mathcal{K}};u,\hbar)$ of $M_{\mathcal{K}}$. As to quantum
invariant $Z(M_{\mathcal{K}};u,\hbar)$, there are also two other
approaches to compute the expansion coefficients $S_n(u)$: the state
integral model introduced in \cite{DGLZ} and the topological
recursion method which is from string theory \cite{DFM}. These two
methods may also be used to study the geometric branch of cases in
ii) which will be further studied in \cite{Zhu}.

\proof[Acknowledgements] This work was finished when the author was
visiting University of California, Los Angeles. The author would
like to thank the China Scholarship Council for the finicial support
for his visit. The author thanks Professor H. Fuji for providing him
partial calculations in their paper \cite{DFM}. The author also
would like to thank Professor Kefeng Liu for his careful reading the
draft of this paper and his valuable suggestions.

$$ \ \ \ \ $$

\end{document}